\documentclass[11pt,twoside]{article}
\usepackage[english]{babel}
\usepackage{amssymb,amsmath}
\usepackage{fancyhea}
\usepackage{amssymb}
\usepackage[width=132mm,height=190mm]{geometry}
\usepackage{enumerate}
\usepackage{amssymb}
\usepackage[english]{babel}
\usepackage[T1]{fontenc}
\usepackage[latin1]{inputenc}

\newtheorem{thm}{Theorem}[section]

\usepackage{amsfonts,latexsym,graphicx}
\usepackage{times}
\usepackage{enumerate}

\newcommand{\F}{\mathbb{F}}
\newcommand{\N}{\mathbb{N}}
\newcommand{\Z}{\mathbb{Z}}
\newcommand{\G}{\mathbb{G}}

\newcommand{\Hom}{\mathrm{hom}}
\newcommand{\fe}{f.\,e.}
\providecommand{\abs}[1]{\lvert#1\rvert}

\begin{document}

\vspace*{2mm}

\noindent {\bf \LARGE New arcs in projective Hjelmslev planes over
\\[3pt] Galois rings}\thispagestyle{fancyplain} \setlength\partopsep
{0pt} \flushbottom
\date{}

\vspace*{4mm}

\noindent {\sc Michael Kiermaier, Axel Kohnert} \\
{\hspace*{12mm} \hfill {\tt
michael.kiermaier,axel.kohnert@uni-bayreuth.de} \\
Mathematical Department, University of Bayreuth, D-95440 Bayreuth,
GERMANY}

\begin{center}
\parbox{11,8cm}{\footnotesize{\bf Abstract.}
It is known that some good linear codes over a finite ring
($R$-linear codes) arise from interesting point constellations in
certain projective geometries. For example, the expurgated
Nordstrom-Robinson code, a nonlinear binary $[14,6,6]$-code which
has higher minimum distance than any linear binary $[14,6]$-code,
can be constructed from a maximal $2$-arc in the projective
Hjelmslev plane over $\Z_4$.

We report on a computer search for maximal arcs in projective
Hjelmslev planes over proper Galois rings of order $\leq 27$. The
used method is to prescribe a group of automorphisms which shrinks
the problem to a computationally feasible size. The resulting
system of Diophantine linear equations is solved by lattice point
enumeration.

We improve many of the known lower bounds on the size of maximal
arcs. Furthermore, the Gray image of one of the constructed arcs
yields a nonlinear quaternary $[504,6,376]$-code. This code has
higher minimal distance than any known $\F_4$-linear
$[504,6]$-code.}
  \end{center}

\setcounter{equation}{0} \setcounter{page}{112}

\vspace*{-2mm}

\section{Galois Rings}

\vspace*{-2mm}

\baselineskip=0.92\normalbaselineskip

For a prime power $q=p^r$ and a natural number
$m\in\N\setminus\{0\}$, the \emph{Galois ring}
$\mbox{GR}(q^{m},p^{m})$ of order $q^m$ and characteristic $p^m$
is defined as $\mathbb{Z}_{p^{m}}[X]/(f)$, where
$f\in\mathbb{Z}_{p^{m}}[X]$ is a monic polynomial of degree $r$
which is irreducible modulo $p$. For different choices of the
polynomial $f$, the resulting Galois rings are isomorphic.

The class of the Galois rings contains the finite fields and the
integers modulo a prime power:

\begin{enumerate}[(i)]

\vspace*{-2mm}

 \item \label{galois:finite_field}
$\mbox{GR}(q,p)\cong\F_q$.

\vspace*{-2mm}

\item \label{galois:integers_modulo}
$\mbox{GR}(p^{m},p^{m})\cong\mathbb{Z}_{p^{m}}$
\end{enumerate}

\vspace*{-2mm}

A Galois ring that is not a finite field will be called
\emph{proper} Galois ring.

The Galois rings are well suited for base rings of linear codes:
Case~(\ref{galois:finite_field}) gives the classical linear codes,
and case~(\ref{galois:integers_modulo}) contains the $\Z_4$-codes.
The smallest Galois ring which is neither a finite field nor a
residue class ring is $\G_{16}:=\mbox{GR}(16,4)$. This ring admits
very good codes, too: One example can be found in
\cite{Hemme-Honold-Landjev-2000}, a new one will be given below.
As a subset of the finite chain rings, we can apply the theory in
\cite{Honold-Landjev-2000-EJC7:R11} to linear codes over Galois
rings, including a generalized Gray isometry
\cite{Greferath-Schmidt-1999-IEEETIT45[7]:2522-2524}.

\vspace*{-2mm}

\section{Arcs in projective Hjelmslev planes}

\lhead{} \rhead{}
\chead[\fancyplain{}{\small\sl\leftmark}]{\fancyplain{}{\small\sl
\leftmark}} \markboth{\hspace{-3,7cm}Fifth  International Workshop
on Optimal Codes and Related Topics\\ June 16-22, 2007, White
Lagoon, Bulgaria  \hfill pp. 112-119}{}

\vspace*{-2mm}

The \emph{projective Hjelmslev plane} $\mbox{PHG}(2,R)$ over a
Galois ring $R=\mbox{GR}(q^m,p^m)$ is defined as follows: The
point set $\mathfrak{P}$ (line set $\mathfrak{L}$) is the set of
the free rank $1$ (rank $2$) submodules of the module $R^3$, and
the incidence is given by set inclusion.

We have $\abs{\mathfrak{P}} = \abs{\mathfrak{L}} = (q^2 + q + 1)
q^{2(m-1)}$. For $m\neq 1$, a projective Hjelmslev plane is not a
classical projective plane, because two different lines may meet
in more than one point. More about projective Hjelmslev geometries
can be found in \cite{Honold-Landjev-2001-DM231:265-278} and the
references cited there.

\pagestyle{myheadings} \markboth{OC2007}{Kiermaier, Kohnert}

 For
$n\in\N$, a set of points $\mathfrak{k}\subseteq\mathfrak{P}$ of
size $n$ is called \emph{projective $(n,u)$-arc}, if some $u$
elements of $\mathfrak{k}$ are collinear, but no $u+1$ elements of
$\mathfrak{k}$ are collinear. If we allow $\mathfrak{k}$ to be a
\emph{multi}set of points in this definition\footnote{Of course we
have to respect multiplicities for counting the number of
collinear points.}, $\mathfrak{k}$ is called
\emph{$(n,u)$-multiarc}.

If $R$ is a finite field, $\mbox{PHG}(2,R)$ is a classical
projective plane. In this case, the arc problem and related
problems are heavily investigated, see
\fe~\cite{Ball_arcs_blocking_sets_1996,BKW_arcs_2005,Hirschfeld_Storme_packing_update_2001}.
To exclude the classical case from the search, we restrict ourself
to proper Galois rings $R$.

\cite{Honold-Landjev-2001-DM231:265-278} contains a table for arcs
over chain rings of composition length $m=2$ and order $\leq 25$.
A few new arcs can be found in \cite{Hemme-Honold-Landjev-2000},
and further improvements for chain rings of composition length
$m=2$ and order $9$ and $25$ are published in
\cite{Boumova-Landjev-2004}. In \cite{Kiermaier-2006} a complete
classification of $(n,u)$-multiarcs was done for small $u$ in
small Hjelmslev geometries over chain rings $\neq \Z_{16}$ of
order $\leq 16$, which again yielded some improvements of the
bounds. The most important results of this search can also be
found in \cite{Honold-Kiermaier-2006}.

\vspace*{-2mm}

\section{Solving Linear Diophantine Equations}

\vspace*{-2mm}

The construction of discrete objects using incidence preserving
group actions is a general approach that works in many cases
\cite{Kerber_finite_group_ed2_1999}. It was first applied in the
70's for the construction of designs
\cite{Kramer_Mesner_ihre_matrix}. Later this method was used for
the construction of $q$-analogs of designs
\cite{Braun_kerber_laue_q_analoga,Braun_q_analoga}, parallelisms
in projective geometries \cite{Braun_point_cyclic_2006}, distance
optimal codes \cite{BraunKohnertWassermann:05} and arcs over
projective planes \cite{BKW_arcs_2005}.

For our problem, we study the \emph{incidence matrix} $M$ of
$\mbox{PHG}(2,R)$: The columns are labeled by the points and the
rows are labeled by the lines. The entry of $M$ indexed by the
line $L\in \mathfrak{L}$ and the point $P\in \mathfrak{P}$ is
defined as
\[
M_{L,P}:=\begin{cases}
  1 & \mbox{if }P\subset L\\
  0 & \mbox{otherwise}
\end{cases}
\]

\vspace*{-2mm}

Using this incidence matrix we can restate the problem of finding
an $(n,u)$-arc as follows:

\begin{thm}
\label{thm:ohne automorphismen} There is a projective $(n,u)$-arc
in $\mbox{PHG}(2,R)$ if and only if there is a $0/1$-solution
$x=(x_{1},\ldots,x_{|\mathfrak{P}|})$ of the following system of
(in)equalities
\[
\begin{array}{ccccc}
(1) & \sum_{i=1}^{\abs{\mathfrak{P}}} x_{i} & = & n\\
(2) & Mx^{T} & \le & \left(\begin{array}{c}
u\\
\vdots\\
u\end{array}\right)\end{array}
\]
and at least one of the lines of the system $(2)$ is an equality.
\end{thm}

This comes from the fact that the entries equal to one in a solution
vector $x$ define the selection of points which go into the arc.
The restriction to a zero-one solution ensures that we get a projective
arc. If one admits non-negative integers for the entries of $x$, the resulting
solutions will be multiarcs.

To solve this system for interesting cases we use lattice point
enumeration based on the \emph{LLL}-algorithm
\cite{wassermann-designs-enumeration-techniques-1998}. But to get
new results we have to solve systems of sizes which are too large
for the solving algorithm. (e.g. $\abs{\mathfrak{P}}=775$ for
$R=\mathbb{Z}_{25}$). To reduce the size of the system we
prescribe automorphisms $\phi\in\mbox{GL}(3,R)$, so we are looking
for solutions (i.e. arcs $\mathfrak{k}$) with the additional
property that
\[
P\in \mathfrak{k}\Rightarrow\phi(P)\in\mathfrak{k}
\]
This means for the matrix $M$ of the system that we can sum up columns
which correspond to points lying in the same orbit. As the defining
incidence property of the matrix $M$ is invariant under the prescribed
automorphism, i.e.
\[
P\subset L\Rightarrow\phi(P)\subset\phi(L)
\]
we get (after the fusion of points which are in the same orbit)
identical rows in the matrix. These are the rows corresponding to
the lines in the orbits, which we get by applying the
automorphisms to the lines. Therefore we can also reduce the
number of rows of the matrix $M$. As the number of orbits is
identical on lines and points, the reduced matrix is again a
square matrix of size $m$ which is the number of orbits. We call
this new matrix $M^{G}$ where $G$ is the group generated by the
prescribed automorphisms. The rows are indexed by the orbits
$\Omega_{1},\ldots,\Omega_{m}$ of the lines, and the columns are
labeled by the orbits $\omega_{1},\ldots,\omega_{m}$ of the
points. An entry of $M^{G}$ is given by
\[
M_{\Omega_{i},\omega_{j}}^{G}:=\abs{\{P\in w_{j}:P\subset L\}}
\]
where $L$ is a representative of $\Omega_{i}$. Now we can restate
the above theorem:

\begin{thm}
There is an $(n,u)$-arc in $\mbox{PHG}(2,R)$ whose
automorphism-group $H$ contains the group $G<\mbox{GL}(3,R)$ as a
subgroup, if and only if, there is a $0/1-$solution
$x=(x_{1},\ldots,x_{m})$ to the following system of
(in)equalities\[
\begin{array}{ccccc}
(1) & \sum|\omega_{i}|x_{i} & = & n\\
(2) & M^{G}x^{T} & \le & \left(\begin{array}{c}
u\\
\vdots\\
u\end{array}\right)\end{array}\]
and at least one of the lines of the system $(2)$ is an equality.
\end{thm}

For computational purposes we transform the system of inequalities
into a system of equations. We solve the following system (we
denote by $-I$ the negative unit matrix, and $(M^{G},-I)$ denotes
the $m\times2m$ block matrix):

\vspace*{-2mm}

\[
\begin{array}{ccccc}
(1) & \sum|\omega_{i}|x_{i} & = & n\\
(2) & (M^{G},-I)(x,y)^{T} & = & \left(\begin{array}{c}
0\\
\vdots\\
0\end{array}\right)\end{array}\]

\vspace*{-2mm}

\noindent The additional variables $y=(y_{1},\ldots,y_{m})$ in a
solution may have values in $\{0,\ldots,u\}.$ From these values we
obtain the intersection numbers between an arc and the lines in
$\mbox{PHG}(2,R)$, so we easily get the secant distribution of an
arc from these values. This is the system of Diophantine equations
we finally solve to get new arcs.

\vspace*{-2.5mm}

\section{An example}

\vspace*{-2mm}

 We constructed a projective $(126,8)$-arc in
$\mbox{PHG}(2,\G_{16})=(\mathfrak{P},\mathfrak{L})$. This plane
has $336$ points and lines, so the incidence matrix $M$ has
$336^2$ entries which would it make difficult to solve the
corresponding system of equations given in theorem \ref{thm:ohne
automorphismen}.

To reduce the size, we take a matrix $A\in\G_{16}^{3\times 3}$
such that its image $\bar{A}\in\F_4^{3\times 3}$ modulo 2
generates a Singer cycle on $\mbox{PG}(2,\F_4)$. So
$\mbox{ord}(A)=21k$ where $k\in\N$. By replacing $A$ by $A^k$ we
can assume $\mbox{ord}(A)=21$. Let $G$ the cyclic group generated
by $A$.

$G$ operates on the point set $\mathfrak{P}$ by
left-multiplication and partitions $\mathfrak{P}$ into $16$ orbits
of length $21$. Each orbit contains exactly one point of each
neighborhood class of $\mbox{PHG}(2,\G_{16})$, so each
neighborhood class contains exactly $6$ points of $\mathfrak{k}$.

The smaller system corresponding to the $16\times 16$ matrix
$M^{G}$ could be solved in a few seconds and we got a solution
$(x,y)$, which necessarily has exactly $6$ entries from the first
part $x$ equal to $1$. The second part $y$ of the solution
contains the secant distribution. The entries of $y$ are either
$0$ or $8$, so this $(126,8)$-arc has only two different
intersection numbers. There is exactly one entry of $y$ equal to
zero, so the arc intersects with $15\cdot 21=315$ lines in $8$
points and with $21$ lines it has no intersection.

There is a further remarkable property of this arc $\mathfrak{k}$:
Each orbit of $G$ on $\mathfrak{P}$ is a maximal $(21,2)$-arc, a
socalled hyperoval. So $\mathfrak{k}$ can be split into $6$
hyperovals. The hyperovals in $\mbox{PHG}(2,\G_{16})$ are unique
up to geometric isomorphism \cite{Kiermaier-2006}. More on their
structure can be found in
\cite{Honold-Landjev-2005-FFA11:292-304}.

\vspace*{-2.5mm}

\section{A nonlinear quaternary $[504,6,376]$-code}

\vspace*{-2mm}

\baselineskip=0.88\normalbaselineskip

We take the homogeneous coordinates of the above $(126,8)$-arc
$\mathfrak{k}$ and put them as columns into a generator matrix of
a $\G_{16}$-linear code $C$. We equip $\G_{16}$ with the
homogeneous weight $w_\Hom:\G_{16}\rightarrow\N$, that is

\vspace*{-2mm}

\[
w_\Hom(a) = \begin{cases} 0 & \mbox{if } a=0\\3 & \mbox{if }a\mbox{ is a unit}\\4 & \mbox{otherwise}\end{cases}
\]

\vspace*{-1mm}

\noindent To calculate the homogeneous weight distribution of the
code $C$, we use the theory and notation in
\cite{Honold-Landjev-2000-EJC7:R11}. Using the information on the
intersection numbers and the neighborhood distribution, we see
that there are only these two $\mathfrak{k}$-types of lines:
\begin{itemize}
\item $21$ lines of $\mathfrak{k}$-type $(a_0,a_1,a_2)=(96,30,0)$
\item $315$ lines of $\mathfrak{k}$-type $(a_0,a_1,a_2)=(96,22,8)$
\end{itemize}
That gives the homogeneous weight enumerator
\begin{multline*}
315\cdot 12\cdot X^{96\cdot 3 + 22\cdot 4} + 21\cdot 12\cdot
X^{96\cdot 3+30\cdot 4} + 21\cdot 3\cdot X^{96\cdot 4} + 1 =
\\ 252 X^{408} + 63 X^{384} + 3780 X^{376} + 1
\end{multline*}
So $C$ is a $\G_{16}$-linear $[126,3,376]$-code. Using the
generalized Gray map $\psi$ as defined in
\cite{Greferath-Schmidt-1999-IEEETIT45[7]:2522-2524}, the code $C$
yields a nonlinear\footnote{We used a computer check to make sure
that $C$ is nonlinear.} distance-invariant quaternary
$[504,6,376]$-code $\psi(C)$.

By applying the one-step Griesmer bound (see
\fe~\cite{Betten-Braun-Fripertinger-Kerber-Kohnert-Wassermann-2006},
page~88), the existence of a linear quaternary $[504,6,376]$-code
would imply the existence of a linear quaternary
$[128,5,94]$-code. According to \cite{Grassl-codetables.de}, no
such code is known.

So $\psi(C)$ clearly is a very good nonlinear code and it might be
better than any linear quaternary code of equal length and size.

\vspace*{-2mm}

\section{New projective arcs}

\vspace*{-2mm}

Table 1 contains the sizes of the arcs we constructed with our
method. All these sizes are at least as big as the previously best
known sizes, with one exception: We only found a projective
$(184,12)$-arc in $\mbox{PHG}(2,\G_{16})$, while the construction
in Example~4.7 in \cite{Honold-Landjev-2001-DM231:265-278} gives a
projective $(186,12)$-arc. This value is indicated by italic font
in the table. If we know that an entry meets some known upper
bound and therefore a corresponding arc is of maximal possible
size, we use bold font. A lower index $*$ denotes an improvement
against the previously known value. We do not use the $*$-symbol
for $\Z_8$, $\Z_{16}$ and $\Z_{27}$, since for these rings only
very few values were published before. An upper index $E$ denotes
an arc we found by extension, this means that we used a program,
which checks whether it is possible to add a further point, which
is not yet in the arc, without violating the defining condition of
maximal $u$ projective collinear points.

%\renewcommand{\topfraction}{1}
%\begin{table}[t]
\begin{center}
\scalebox{0.81} {%\footnotesize
\begin{tabular}{|c||c|r@{}l|r@{}l|r@{}l|r@{}l|r@{}l|}
\hline
$_{u}\backslash^{R}$& $\mathbb{Z}_{8}$& \multicolumn{2}{|c|}{$\mathbb{Z}_{9}$}&
\multicolumn{2}{|c|}{$\mathbb{G}_{16}$}& \multicolumn{2}{|c|}{$\mathbb{Z}_{16}$}& \multicolumn{2}{|c|}{$\mathbb{Z}_{25}$}&
\multicolumn{2}{|c|}{$\mathbb{Z}_{27}$}\tabularnewline
\hline
\hline
$\abs{\mathfrak{P}}$& $112$& \multicolumn{2}{|c|}{$117$} & \multicolumn{2}{|c|}{$336$} & \multicolumn{2}{|c|}{$448$} & \multicolumn{2}{|c|}{$775$} &
\multicolumn{2}{|c|}{$1053$}\tabularnewline
\hline
\hline
$2$& $\mathbf{10}$& $\mathbf{9}$& & $\mathbf{21}$& & $\mathbf{16}$& & $20$& & $21$ & \tabularnewline
\hline
$3$& $\mathbf{21}$& $\mathbf{19}$& & $27$& & $28$& & $34$& ${}_{*}$& $39$ & \tabularnewline
\hline
$4$& $28$& $\mathbf{30}$ & ${}_{*}^{E}$ & $\mathbf{52}$ & ${}_{*}$& $49$& & $60$ & ${}_{*}$& $65$ & \tabularnewline
\hline
$5$& $37$ & $39$ & ${}_{*}$& $68$ & ${}_{*}^{E}$& $61$& & $79$ & ${}_{*}^{E}$ & $91$ &\tabularnewline
\hline
$6$ & $48$ & $48$ & ${}_{*}$& $\mathbf{84}$ & ${}_{*}$& $84$& & $105$ & ${}_{*}^{E}$ & $117$ & \tabularnewline
\hline
$7$ & $56$ & $60$ & ${}_{*}^{E}$& $93$ & ${}_{*}$& $100$& & $124$ & ${}_{*}$ & $130$ & \tabularnewline
\hline
$8$& $68$ & $69$ & ${}_{*}$ & $\mathbf{126}$ & ${}_{*}$ & $120$& & $150$ & ${}_{*}$ & $172$ & ${}^{E}$\tabularnewline
\hline
$9$ & $79$ & $\mathbf{81}$& & $140$ & ${}_{*}$ & $136$ & & $175$ & ${}_{*}$& $193$ & \tabularnewline
\hline
$10$& $88$& $\mathbf{93}$ & & $152$ & ${}_{*}$& $156$& & $199$ & ${}_{*}$& $234$ & \tabularnewline
\hline
$11$& $100$ & $\mathbf{105}$& & $164$ & ${}_{*}^{E}$ & $172$ & & $223$ & ${}_{*}$& $240$ & ${}^{E}$\tabularnewline
\hline
$12$& $\mathbf{112}$& $\mathbf{117}$ & & $\mathit{184}$ & & $208$ & & $256$ & ${}_{*}^{E}$& $288$ & \tabularnewline
\hline
$13$& & & & $200$ & ${}_{*}$& $212$ & ${}^{E}$& $310$ & ${}_{*}$& $302$ & ${}^{E}$\tabularnewline
\hline
$14$& & & & $216$ & ${}_{*}$& $232$& &  $311$ & ${}_{*}^{E}$& $351$ & \tabularnewline
\hline
$15$& & & & $236$ & ${}_{*}$& $252$& & $328$ & ${}_{*}$& $369$ & \tabularnewline
\hline
$16$& & & & $\mathbf{256}$& & $276$& & $355$ & ${}_{*}^{E}$& $390$ & \tabularnewline
\hline
$17$& & & & $\mathbf{276}$& & $292$& & $385$ & ${}_{*}$& $417$ & ${}^{E}$\tabularnewline
\hline
$18$& & & & $\mathbf{296}$& & $312$& & $417$ & ${}_{*}$& $468$ & \tabularnewline
\hline
$19$& & & & $\mathbf{316}$& & $331$ & ${}^{E}$& $465$ & ${}_{*}$& $505$ & \tabularnewline
\hline
$20$& & & & $\mathbf{336}$& & $360$& & $480$ & ${}_{*}$& $520$ & ${}^{E}$\tabularnewline
\hline
$21$& & & & & & $376$& & $510$& & $567$ & \tabularnewline
\hline
$22$& & & & & & $400$& & $534$ & ${}_{*}^{E}$& $595$ & ${}^{E}$\tabularnewline
\hline
$23$& & & & & & $424$& & $565$ & ${}_{*}$ & $617$ & ${}^{E}$\tabularnewline
\hline
$24$& & & & & & $\mathbf{448}$& & $592$ & ${}_{*}$& $657$ & \tabularnewline
\hline
$25$& & & & & & & & $\mathbf{625}$& & $676$ & ${}^{E}$\tabularnewline
\hline
$26$& & & & & & & & $\mathbf{655}$& & $702$ & \tabularnewline
\hline
$27$& & & & & & & & $\mathbf{685}$& & $747$ & \tabularnewline
\hline
$28$& & & & & & & & $\mathbf{715}$& & $783$ & \tabularnewline
\hline
$29$& & & & & & & & $\mathbf{745}$& & $819$ & ${}^{E}$\tabularnewline
\hline
$30$& & & & & & & & $\mathbf{775}$& & $840$ & \tabularnewline
\hline
$31$& & & & & & & & & & $873$ &\tabularnewline
\hline
$32$& & & & & & & & & & $909$ &\tabularnewline
\hline
$33$& & & & & & & & & & $945$ &\tabularnewline
\hline
$34$& & & & & & & & & & $981$ &\tabularnewline
\hline
$35$& & & & & & & & & & $\mathbf{1017}$ &\tabularnewline
\hline
$36$& & & & & & & & & & $\mathbf{1053}$ &\tabularnewline
\hline
\end{tabular}}
\label{result_table}
\\[12pt]
Table 1. Sizes of the constructed projective $(n,u)$-arcs
\end{center}
%\end{table}

%\vspace*{-2mm}

All the values in the table are for projective arcs. However, we
know a non-pro\-jective multiarc which is bigger than the best
known projective arc. This is a $(155,8)$-multiarc in
$\mbox{PHG}(2,\Z_{25})$, it can be constructed in a similar way as
the arc in the above example: The action of a lifted Singer cycle
splits the point set of $\mbox{PHG}(2,\Z_{25})$ into $25$ orbits
of length $31$. It is possible to select $4$ of these orbits, one
of them twice, such that they together give the
$(155,8)$-multiarc.

More information on the arcs (secant distribution, used group of
automorphisms) can be found on the home pages of the authors.

\end{document}